\documentclass{article}
\usepackage{amsmath}
\usepackage{amssymb}
\usepackage{color}

\usepackage{latexsym}
\usepackage[english]{babel}
\usepackage[dvips]{graphicx}
\usepackage{verbatim}
\newcommand{\fintersect}{forwardly intersect}
\newcommand{\bad}{unbounded}

\newcommand{\Supp}{\mathit{Supp}}
\newcommand{\dbt}{(X+Y)/2}
\newcommand{\Xt}{X'}
\newcommand{\Yt}{Y'}
\newcommand{\TT}{V_p}

\newcommand{\bbibitem}{\bibitem}

\newcommand{\llabel}[1]{{\label{#1}}}

\renewcommand{\r}[1]{(\ref{#1})}

\newcommand{\ex}[1]{} 

\font\tenmsb=msbm10
\font\sevenmsb=msbm7
\font\fivemsb=msbm5
\newfam\msbfam
\textfont\msbfam=\tenmsb
\scriptfont\msbfam=\sevenmsb
\scriptscriptfont\msbfam=\fivemsb

\newcommand{\bi}{\begin{itemize}}
\newcommand{\ei}{\end{itemize}}
\newcommand{\bd}{\begin{description}}
\newcommand{\ed}{\end{description}}
\renewcommand{\i}{\item}

\newcommand{\bqn}{\begin{eqnarray}}
\newcommand{\eqn}{\end{eqnarray}}
\newcommand{\eqnn}{\nonumber\end{eqnarray}}
\newcommand{\eqnl}[1]{\llabel{#1}\end{eqnarray}}

\newcommand{\nn}{\nonumber}
\newcommand{\ba}[1]{\begin{array}{#1}}
\newcommand{\ea}{\end{array}}

\newcommand{\R}
{\mathbb{R}}

\newcommand{\N}{\mathbb{N}}

\newcommand{\fine}{\end{document}}

\def \trait (#1) (#2) (#3){\vrule width #1pt height #2pt depth #3pt}
\def \qed{\hfill
        \trait (0.1) (6) (0)
        \trait (6) (0.1) (0)
        \kern-6pt   
        \trait (6) (6) (-5.9)
        \trait (0.1) (6) (0)
\medskip}

\newtheorem{Theorem}{\bf Theorem}
\newtheorem{ml}[Theorem]{\bf Lemma}
\newtheorem{mo}{\bf \underline{{\sl Observation}}}
\newtheorem{mcc}[Theorem]{\bf Corollary}
\newtheorem{Definition}[Theorem]{\bf Definition}
\newtheorem{mpr}[Theorem]{\bf Proposition}
\newtheorem{mproperty}[Theorem]{\bf Property}
\newtheorem{corollary}[Theorem]{Corollary}
\newtheorem{mrem}{\bf \underline{{\sl Remark}}}
\newtheorem{rmk}[Theorem]{Remark}

\newcommand{\bt}{\begin{Theorem}}
\newcommand{\et}{\end{Theorem}}
\newcommand{\bl}{\begin{ml}}
\newcommand{\el}{\end{ml}}
\newcommand{\bo}{\noindent\begin{mo}\rm}
\newcommand{\eo}{\end{mo}}
\newcommand{\bp}{\begin{mpr}}
\newcommand{\ep}{\end{mpr}}
\newcommand{\bc}{\begin{mcc}}
\newcommand{\ec}{\end{mcc}}
\newcommand{\bdeff}{\begin{Definition}}
\newcommand{\edeff}{\end{Definition}}
\newcommand{\bproperty}{\begin{mproperty}}
\newcommand{\eproperty}{\end{mproperty}}
\newcommand{\brem}{\begin{mrem}\rm}
\newcommand{\erem}{\end{mrem}}

\newcommand{\bfi}{\begin{figure}}
\newcommand{\efi}{\end{figure}}

\newcommand{\capt}{
\refstepcounter{figure}
\begin{center}Figure~\thefigure\end{center}
}

\newcommand{\lb}{\lambda}

\newcommand{\al}{\alpha}
\newcommand{\eps}{\varepsilon}
\newcommand{\de}{\delta}

\newcommand{\Om}{\Omega}

\newcommand{\con}{{\cal C}}
\newcommand{\E}{{\cal E}}

\newcommand{\A}{{\cal A}}

\newcommand{\U}{{\cal U}}
\newcommand{\T}{{\cal T}}

\newcommand{\qq}{{\cal Z}}

\newcommand{\D}{{\cal D}}

\newcommand{\proof}{{\bf Proof. }}

\newcommand{\brs}{\begin{eqnarray*}}
\newcommand{\ers}{\end{eqnarray*}}
\newcommand{\br}{\begin{eqnarray}}
\newcommand{\er}{\end{eqnarray}}

\newcommand{\rw}{\rightarrow}
\newcommand{\Ga}{\Gamma}
\newcommand{\ga}{\gamma}
\newcommand{\lp}{\left(}
\newcommand{\rp}{\right)}
\newcommand{\mt}{\mapsto}
\def\EOP{\ \hfill \rule{0.5em}{0.5em} }
\newcommand{\be}{\begin{equation}}
\newcommand{\ee}{\end{equation}}
\newcommand{\n}{\noindent}
\newcommand{\coco}{{\cal G}}

\newcommand{\C}{\mathbb{C}}

\newcommand{\B}{{\cal B}}

\newcommand{\Bbq}{\overline{\B(q)}}
\newcommand{\Dbq}{\overline{{\cal D}(q)}}

\begin{document}
\begin{center} \noindent
{\LARGE{\sl{\bf Stability of Planar Nonlinear Switched Systems}}}
\end{center}

\vskip 1cm
\begin{center}
Ugo Boscain,

{\footnotesize SISSA-ISAS,
Via Beirut 2-4, 34014 Trieste, Italy. E-mail: {\tt boscain@sissa.it}}

Gr\'egoire Charlot

{\footnotesize Universit\'e Montpellier II, Math\'ematiques, CC051, 34095 Montpellier Cedex 5, 
France. E-mail: {\tt cha@math.univ-montp2.fr }}

Mario Sigalotti

{\footnotesize  INRIA, 2004 route des lucioles, 06902 Sophia Antipolis, France. E-mail: {\tt Mario.Sigalotti@sophia.inria.fr}}

\end{center}

\vspace{1cm}

\begin{quotation}\noindent  {\bf\em Abstract --- 
We consider the time-dependent nonlinear system
$\dot q(t)=u(t)X(q(t))+(1-u(t))Y(q(t))$, where $q\in\R^2$, $X$ and
$Y$ are two 
smooth vector fields, globally  asymptotically
stable at the origin and
  $u:[0,\infty)\to\{0,1\}$ is an arbitrary measurable function.
Analysing the topology of the set where $X$ and $Y$ are parallel,
we give some sufficient and some necessary conditions for global
asymptotic stability, uniform with respect to $u(.)$. Such conditions
can be verified without any integration or construction of a Lyapunov
function, and they are
robust under small perturbations of the vector fields.
}\end{quotation}


\section{Introduction}
A  {\it switched system} is a family of continuous-time dynamical
systems endowed with a 
rule that determines, at every time, which dynamical
system is responsible for the time evolution. More precisely let $\{f_u\,|\;u\in U\}$
be a
(possibly infinite) set of smooth vector fields  on
a manifold $M$, and consider, as $u$ varies in $U$, the family of dynamical systems
\bqn\label{-1}
\dot q=f_u(q)\,,~~~~q\in M\,.
\eqn
A non-autonomous dynamical system is obtained by assigning 
a so-called {\it switching
function} $u(.):[0,\infty)\to U$. 

In this paper, the switching function models the behavior of a parameter 
which cannot be predicted a priori. It
represents some phenomena (e.g., a
disturbance) that it is not possible to control or include in the
dynamical system model.

A typical problem related to switched systems is to obtain, out of a 
property which is shared by all the autonomous dynamical systems governed
by the vector fields  $f_u$, some, maybe weaker, property for the time-dependent system
associated with an arbitrary switching function $u(.)$. For a discussion
on various issues related to switched systems we refer the reader to
\cite{liberzon-book,survey}.

In this paper, we consider a two-dimensional
nonlinear switched system of the type
\bqn
\dot q=u\,X(q)+(1-u)\,Y(q)\,,~~~~q\in\R^2\,,~~~~u\in\{0,1\}\,,
\label{system}\label{nonlinear}
\eqn
where the two vector fields 
$X$ and $Y$ are smooth (say, $\con^\infty$) on $\R^2$.
In order to define proper non-autonomous systems, we require
the switching functions to be measurable.

Assume that $X(0)=Y(0)=0$ and that the two dynamical
systems $\dot q=X(q)$ and $\dot q=Y(q)$ are \underline{globally 
asymptotically
stable} at the origin.
Our main aim is to study under which conditions on $X$
and $Y$ the origin is globally asymptotically stable
for the system \r{nonlinear}, uniformly with respect to the switching functions
(GUAS for short). For the precise formulation of this and other stability properties, 
see Definition \ref{stabilities}. 

In order to study the stability of \r{system} it is natural to consider 
its convexification, i.e., the case in which $u$ varies in the whole interval $[0,1]$. It  turns 
out that the stability properties of the two systems are equivalent (see Section 
\ref{s-convex}).\\\\
The linear version of the system introduced above, 
namely,  
\bqn\label{1}
\dot q=u\,A\,q+(1-u)\,B\,q\,,~~~~q\in\R^2\,,~~~~u\in\{0,1\}\,,
\eqn
where  the $2\times2$ real matrices
$A$ and $B$ have eigenvalues with 
strictly negative real part, 
 was studied in \cite{SIAM} (see also \cite{lyapunov-paolo}). More
precisely, the results in \cite{SIAM} establish
a necessary and sufficient condition for GUAS 
in terms of three relevant parameters, two depending 
on the eigenvalues of $A$ and $B$ respectively, and the third one 
(namely, the cross ratio of the four eigenvectors of $A$ and $B$ in
the projective line $\C P^1$) accounting for the interrelations
among the two systems.
The precise necessary and sufficient condition ensuring GUAS 
of 
\r{1} is quite
technical and can be found in \cite{SIAM} (see also \cite{lyapunov-paolo}).
Notice that, in the linear case, 
GUAS
is equivalent to
the more often quoted GUES property, i.e., global exponential 
stability, uniform with respect to the switching rule (see, for example, 
\cite{angeli} and references therein). 
For related results on linear switched systems, see \cite{agr,blanchini,daw,hesp,lyapunov-paolo}.

\bigskip

For nonlinear systems, the problem of characterizing GUAS completely, 
without assuming the explicit knowledge of the 
integral curves of $X$ and $Y$, 
is hopeless. 

The problem, however, admits some partial solution. 
The purpose of this paper is to provide some sufficient and some necessary 
conditions for stability which are robust (with respect to 
small perturbations of the vector fields) and easily verifiable, directly on 
the vector fields $X$ and $Y$, without requiring any integration or 
construction of a Lyapunov function. 

Denote by $\qq$ the set on which $X$ and $Y$ are parallel.
One of our main results is that, if $\qq$ reduces to the singleton 
$\{0\}$, then \r{nonlinear} is GUAS (Theorem \ref{t-mai-paralleli}). 
The proofs works by showing that an admissible trajectory starting from a point $p\in\R^2$ is forced 
to stay in a compact region bounded by the integral curves of $X$ and $Y$ from $p$. The fact that $X$ and $Y$ are linearly independent  outside the origin plays as a sort of drift
which guarantees that the only possible accumulation point of an admissible trajectory is the origin.

When $\qq$ is just compact, we prove that \r{system} is at least
bounded (see Theorem \ref{t-solo-compatte}). Roughly speaking, this means that its trajectories do not 
escape to infinity. The idea of the proof is that, if we modify $X$ and $Y$ only in a compact
region of the plane, then the boundedness properties of the system are left unchanged. 
Taking advantage of the result obtained in Theorem~\ref{t-mai-paralleli}, we manage to 
prove the boundedness of \r{system} by reducing, using compact perturbations, $\qq$ to $\{0\}$, 
while preserving the global asymptotic stability of $X$ and $Y$.

Other conditions can be formulated taking into account  
the relative position of $X$ and $Y$ along $\qq$. Assume that 
$\qq\setminus\{0\}$ contains at least one point $q_0$. Since  
both $X(q_0)$ and $Y(q_0)$ are different from zero, 
the property of pointing in the same or in the opposite versus can be
stated unambiguously.
If  $X(q_0)$ and $Y(q_0)$
have opposite versus, then there exists a switching function, for the 
convexified system, whose output
is the constant trajectory which stays in $q_0$. As a consequence, the 
system \r{system} is not GUAS.

Additional results can be obtained under the assumption that the
pair of vector fields $(X,Y)$ is generic. 
(For the notion of genericity appropriate to our aims, see Section \ref{s-basic-facts}.) 
In particular, the genericity assumption can be used to guarantee 
that 
$\qq\setminus\{0\}$ is  an embedded one-dimensional 
submanifold of the plane.
Clearly, $\qq$ needs not to be connected. 
If the connected component of $\qq$ containing the origin reduces to $\{0\}$ 
and on all other components $X$ and $Y$ point in the same versus, transversally 
to $\qq$, then 
\r{system} is GUAS. This result is formulated in Theorem \ref{t-mai-tangent}, 
which follows the pattern of proof of Theorem~\ref{t-mai-paralleli}.

Conversely, Theorem \ref{t-inv-non-comp} states that, if one 
connected component of $\qq\setminus\{0\}$ is unbounded and such that 
$X$ and $Y$ have opposite versus on it, 
then  \r{system}  admits a trajectory going to infinity.  
Intuitively, this happens because 
the orientation of $(X(p),Y(p))$ changes while $p$ crosses $\qq\setminus\{0\}$. If $X(p)$
is not tangent to $\qq$ at $p$ and $X(p)$ points in the opposite direction with respect to $Y(p)$, 
then one can embed $\qq$, locally near $p$, in a foliation 
made of admissible trajectories of \r{system}, whose 
running direction 
is reversed while crossing $\qq$ (see Figure~\ref{foli}).
\bfi[h]
\begin{center}
\input{foli.pstex_t}
\caption{A local foliation embedding $\qq$}\label{foli}
\end{center}
\efi
Since, generically, the points where $X$ is tangent to $\qq$ are isolated, it turns out that 
there exists
an admissible trajectory which tracks globally the unbounded connected component of $\qq\setminus\{0\}$
on which $X$ and $Y$ have opposite versus.

\bigskip
The paper is organized as follows.
In Section \ref{s-basic-facts}, we recall the main definitions of stability in which we are 
interested, we introduce the convexified system, and we describe
the topological structure of the  set $\qq$. 
The main results are stated in Section \ref{s-main-results}, where their robustness is also discussed.
The proofs are given in Sections~\ref{mp}, \ref{mt}, \ref{sc}, and \ref{inc}.

\section{Basic definitions and facts}\label{s-basic-facts}
\subsection{Definitions of stability}
Fix $n,m\in \N$ and consider the switched system
\bqn\label{ss-system}
\dot q=f_u(q)\,,~~~~q\in \R^n\,,~~~~u\in U\subset\R^m\,,
\eqn
where $U$ is a measurable subset of  $\R^m$ and $(q,u)\mapsto f_u(q)$ is 
the restriction on  $\R^n\times U$ of a $\con^\infty$ function from $\R^n\times\R^m$ to 
$\R^n$. Assume that $f_u(0)=0$ for every $u\in U$. \\
\noindent
For every $\de>0$, denote by $B_\de\subset \R^n$ the ball of radius $\de$, centered at the 
origin.
Set
$$\U=\{u:[0,\infty)\rw U\,|\;u(.)\ \mbox{measurable}\}\,.$$

For every $u(.)$ in $\U$ and every $p\in \R^n$, 
denote by $t\mt \ga(p,u(.),t)$ the solution of 
\r{ss-system} such that $\ga(p,u(.),0)=p$. 
Notice that, in general, $t\mt \ga(p,u(.),t)$ needs not to be 
defined for every $t\geq 0$, since the non-autonomous vector field 
$f_{u(t)}$ may not be complete. 
Denote by $\T(p,u(.))$ the maximal element of $(0,+\infty]$ such that 
$t\mt \ga(p,u(.),t)$ is defined on 
$[0,\T(p,u(.)))$, and let 
$$\Supp(\ga(p,u(.),.))=\ga(p,u(.),[0,\T(p,u(.))))\,.$$
If $\Supp(\ga(p,u(.),.))$ is bounded, then $\T(p,u(.))=
+\infty$.

Given 
$p\in \R^n$, the {\it accessible set from $p$}, 
denoted by ${\A(p)}$, is 
defined as
$$\A(p)=\cup_{u(.)\in\U}\Supp(\ga(p,u(.),.))\,.$$ 

\n Several notions of stability for the switched system \r{ss-system} can be 
introduced.
\begin{Definition}\label{stabilities}
We say that \r{ss-system} is
\bi
\i {\bf \bad} 
   if there exist $p\in\R^n$ and $u(.)\in\U$ such that
 $\ga(p,u(.),t)$ goes to infinity as $t$ tends to $\T(p,u(.))$;

\i {\bf bounded} 
   if, for every $K_1\subset \R^n$ compact, there exists $K_2\subset \R^n$ 
compact such that $\ga(p,u(.),t)\in K_2$ for every $u(.)\in\U$, $t\geq 0$ and $p\in K_1$;


\i {\bf uniformly stable}  at the origin 
if, for every 
$\de>0$, there exists $\eps>0$ such that $\A(p)\subset B_\de$ for every 
$p\in B_\eps$;

\i {\bf locally attractive}  at the origin if there exists $\de>0$
such that, for every $u(.)\in\U$ and every $p\in B_\de$,  $\ga(p,u(.),t)$ 
converges to the origin 
as $t$ goes to infinity;

\i {\bf globally attractive}  at the origin if, 
for every $u(.)\in\U$ and every $p\in\R^n$, $\ga(p,u(.),t)$ 
converges to the origin 
as $t$ goes to infinity; 

\i {\bf globally uniformly attractive}  at the origin if, 
for every $\de_1, \de_2>0$, there exists $T>0$ such that  $\ga(p,u(.),T)\in B_{\de_1}$ for every 
$u(.)\in\U$ and every $p\in B_{\de_2}$;

\i {\bf globally uniformly stable (GUS)}  at the origin 
   if it is bounded and uniformly stable at the origin;

\i {\bf locally asymptotically stable (LAS)}  at the origin if it is uniformly stable and locally attractive at the origin;

\i {\bf globally asymptotically stable (GAS)}  at the origin if it is 
uniformly stable and globally attractive at the origin;

\i {\bf globally uniformly asymptotically stable (GUAS)}  at the origin if it is uniformly stable and globally uniformly  attractive at the origin.
\ei
\end{Definition}

It has been showed by Angeli, Ingalls, Sontag, and Wang \cite{AISW} that, 
when $U$ is compact, the notions of  GAS and GUAS are equivalent. This is 
the case for system \r{system}. Moreover, it is well known that, in 
the case in which all the vector fields $f_u$ are linear,  local and 
global properties are equivalent.

\subsection{The convexified system}
\label{s-convex}
In this paper, we focus on the planar switched system 
\be\llabel{nons-system}
\dot q=u\,X(q)+(1-u)\,Y(q)\,,~~~~q\in\R^2\,,~~~~u\in\{0,1\}\,, 
\ee
where $X$ and $Y$ denote two 
vector fields on $\R^2$, of class $\con^\infty$, such
that $X(0)=Y(0)=0$. 
We assume moreover that $X$ and $Y$ are \underline{globally asymptotically 
stable} at the origin. 
Notice, in particular, that $X$ and $Y$ are forward complete.

\bigskip

A classical tool in stability analysis is the convexification of the set of admissible velocities. 
Such transformation does not change the closure of the accessible sets. Moreover,  it was proved in 
\cite{12} (see also \cite[Proposition 7.2]{AISW}) that,
for every $p'\in\R^2$, every switching function $u':[0,\infty)\to[0,1]$,
and every positive continuous function $r$ defined on $[0,\T(p',u'(.)))$, there exist
$u(.)\in \U$
and $p\in\R^2$ such that
\[
\|\ga(p,u(.),t)-\ga(p',u'(.),t)\|\leq r(t)
\]
for every $t\in[0,\T(p',u'(.)))$.
As a consequence  each of the notions
introduced in Definition~\ref{stabilities}
holds for (\ref{nons-system}) if and only if it holds for the same system
where
$U=\{0,1\}$ is
replaced by $[0,1]$.

In the following, to simplify proofs, we deal with the convexified system
\be\llabel{s-system}
\dot q=u\,X(q)+(1-u)\,Y(q)\,,~~~~q\in\R^2\,,~~~~u\in[0,1]\,. 
\ee

\smallskip
\noindent
{\bf Notations.} 
When $u(.)$ is constantly equal 
to zero (respectively, one), we write 
$\ga_Y(p,t)$ (respectively, $\ga_X(p,t)$) for $\ga(p,u(.),t)$. 
Given $p,p'\in\R^2$ and $u(.),u'(.)$ in $\U$, we say that $\ga(p,u(.),.)$ 
and $\ga(p',u'(.),.)$ {\it \fintersect} if 
$\Supp(\ga(p,u(.),.))$ and 
$\Supp(\ga(p',u'(.),.))$ have nonempty intersection.

\subsection{The collinearity set of $X$ and $Y$}
\n A key object in order to detect stability properties 
of \r{s-system} turns out to be the set $\qq$ on which $X$ and $Y$ are 
parallel. We have that $\qq=Q^{-1}(0)$, where  
\bqn\label{decca}
Q(p)=\mathit{det}(X(p),Y(p))\,,\ \ \ \ p\in\R^2\,.
\eqn

In \cite{SIAM}, the stability of the linear switched system \r{1} was studied by  
associating with every point of $\R^2$ 
a suitably defined ``worst'' trajectory
passing through it, whose construction was based upon $\qq$.
The global asymptotic stability of the linear switched system \r{1} 
was then proved to be equivalent to the convergence to the origin 
of every such worst trajectory.
We recall that in the linear case, 
excepted for some degenerate
situations, $\qq$ is either equal to $\{0\}$  or 
is made of  
two straight lines passing through
the origin. 

In the nonlinear case, the
situation is more complex.
Let us represent 
$\qq$ as
\bqn\label{dec}
\qq=\{0\}\cup\bigcup_{\Gamma\in \coco}\Gamma\,, 
\eqn
where $\coco$ is the set of all connected components of $\qq\setminus\{0\}$. Notice that $\coco$ needs not, in general, to be 
countable.
With a slight abuse of notation, we will refer to the elements of 
$\coco$ as to the {\it components} of $\qq$.
\bdeff
Let $\Gamma$ be a component of $\qq$ and fix $p\in\Gamma$. We say that 
$\Gamma$ is {\it direct} (respectively, {\it inverse}) if $X(p)$ and $Y(p)$ 
have the same (respectively, opposite) direction. 
\edeff
\begin{rmk}
The definition is independent  of the choice of $p$, since neither $X$ nor $Y$ 
vanish along $\Gamma$.
\end{rmk}
An example of how $\qq$ can look like is represented in Figure~\ref{qq}.

\bfi[h]
\begin{center}
\input{qq.pstex_t}
\caption{The set  $\qq$}\label{qq}
\end{center}
\efi

\bigskip

Some of the results of this paper are obtained assuming that the set $\qq$ has suitable 
regularity properties, which are generic in the sense defined below. 

A base for the {\it Withney topology} on $\con^\infty(\R^2,\R^2)$ 
(the set of 
smooth vector fields on $\R^2$) can be defined, using the multi-index notation, as the 
family of sets of the type
\[
{\cal V}(k,f,r)=\left\{g\in 
\con^\infty(\R^2,\R^2)\,\left|\;\left\|\frac{\partial^{|I|}(f-g)}{\partial 
x^I}(x)\right\|<r(x), \forall x\in\R^2,|I|\leq k\right.\right\}\,,\]
where $k$ is a nonnegative integer, $f$ belongs to  $\con^\infty(\R^2,\R^2)$, and $r$ is 
a positive continuous function defined on $\R^2$.
Denote by GAS$(\R^2)$ the set of smooth vector fields on $\R^2$ which are globally 
asymptotically stable at the origin, and endow it with the topology induced by Withney's one. 
A {\it generic property} for (\ref{s-system}) is a property which holds for 
an open dense subset of GAS$(\R^2)\times$GAS$(\R^2)$, endowed with the product topology of GAS$(\R^2)$.

\bl \label{genio}
For a generic pair of vector fields $(X,Y)$, 
$\qq\setminus\{0\}$ is an embedded one-dimensional submanifold of $\R^2$. 
Moreover, $Q(p)$ changes sign while $p$ crosses $\qq\setminus\{0\}$.
\el
The lemma is a standard result in genericity theory. It follows from the fact that the condition
\bd
\i[
(G1)] If $p\neq0$ and $Q(p)=0$, then $\nabla Q(p)\neq0$,
\ed
is generic (see, for instance, \cite{generiko}). When $\qq\setminus\{0\}$ is a manifold, we say 
that $p\in \qq\setminus\{0\}$ is a \underline{tangency point} if $X(p)$ is tangent to $\qq$. Under 
condition {\bf (G1)}, $p\in \qq\setminus\{0\}$ is a tangency point if and only if  $\nabla Q(p)$ and $X(p)$ (equivalently, $Y(p)$) are 
orthogonal.

Some of our results are obtained under additional generic conditions. 
One of these, namely,
\bd
\i[
(G2)] The Hessian matrix of $Q$ at the origin is non-degenerate,
\ed
ensures that $\qq$, in a neighborhood of the origin, is given either by $\{0\}$ or by the union 
of two transversal one-dimensional manifolds intersecting at the origin.

Under the generic conditions {\bf (G1)} and  {\bf (G2)}, the connected component of 
$\qq$ containing the origin looks like one of Figure~\ref{f-origin}.

A third generic condition which we will sometimes assume to hold is 
\bd
\i[
(G3)] If $p\neq0$, $Q(p)=0$, and $\nabla Q(p)$ is orthogonal to $X(p)$, 
then the second derivative of $Q$ at $p$ along $X$ (equivalently, $Y$) is different from zero,
\ed
which, together with {\bf (G1)}, guarantees that the tangency points on $\qq$ are isolated.

\bfi[h]
\begin{center}
\begin{picture}(0,0)%
\includegraphics{f-origin.pstex}%
\end{picture}%
\setlength{\unitlength}{2131sp}%
\begingroup\makeatletter\ifx\SetFigFont\undefined%
\gdef\SetFigFont#1#2#3#4#5{%
  \reset@font\fontsize{#1}{#2pt}%
  \fontfamily{#3}\fontseries{#4}\fontshape{#5}%
  \selectfont}%
\fi\endgroup%
\begin{picture}(10329,3519)(1639,-4798)
\end{picture}

\caption{The connected component of $\qq$ containing the origin}\label{f-origin}
\end{center}
\efi

\section{Statement of the results}\label{s-main-results}

We organize our results in sufficient 
and necessary conditions with respect to the stability properties. 

Notice that all such conditions are easily verified 
without any integration or construction of a Lyapunov
function. Moreover, they are  robust under small perturbations of the vector fields, as explained 
in Section \ref{s-robustness}. Let us recall that 
$X$ and $Y$ are assumed to be globally asymptotically stable at the origin and that 
all the results given 
below,
although stated for the case $u\in[0,1]$, are also valid 
for the system where $u$ varies in $\{0,1\}$.


Before stating our main theorems, observe that classical results on 
linearization imply the following.
\bp
\label{nec-cond}  
Assume that the eigenvalues of $A=\nabla X|_{p=0}$ and $B=\nabla Y|_{p=0}$ have 
strictly negative real part.
Then \r{s-system} is  LAS if and only if \r{1} is GUAS. 
\ep
\subsection{Sufficient conditions}
The following theorem gives a simple sufficient condition for GUAS, 
which generalizes the analogous one already known 
for the linear system \r{1} (see \cite{SIAM,lyapunov-paolo}).
\bt 
\label{t-mai-paralleli}
Assume that 
$\qq=\{0\}$. Then the switched system \r{s-system} is GUAS 
at the origin.
\et
Under the generic assumptions 
{\bf (G1)} and  {\bf (G2)}, Theorem \ref{t-mai-paralleli} can be generalized as follows. 
\bt\label{t-mai-tangent}
Assume that the  generic conditions {\bf (G1)} and {\bf (G2)} hold. 
Assume, moreover, that  the origin is isolated in $\qq$ and that there is no  tangency point in 
$\qq\setminus\{0\}$. Then the 
switched system \r{s-system} is GUAS.
\et
When $\qq$ is bounded, although different from $\{0\}$, 
some weaker version of Theorem \ref{t-mai-paralleli} still holds.
\bt
\label{t-solo-compatte}
Assume that $\qq$ is compact. Then the switched system \r{s-system} is bounded.
\et
As a direct consequence of Proposition 
\ref{nec-cond} and Theorem \ref{t-solo-compatte}, we have the following 
sufficient condition for GUS.
\begin{corollary}\llabel{c1}
Let $\qq$ be compact, and the linearized 
switched system be non-degenerate and GUAS. 
Then 
the switched system \r{s-system} is GUS. 
\end{corollary}
\subsection{Necessary conditions}
The following proposition expresses  the straightforward remark that
the inverse  components of $\qq$  constitute obstructions to the 
stability of \r{s-system}. 
The reason is clear: if $\Gamma$ is inverse and $p$ belongs to $\Gamma$, 
then a constant switching function $u(.)$ exists 
such  that $\ga(p,u(.),t)=p$ for every $t\geq 0$.
\bp 
\label{p-mai-inverse}
If $\qq$ has an inverse component, then the switched system \r{s-system} is not globally attractive.
\ep
The following theorem gives a necessary condition for boundedness, under generic conditions.
\bt
\label{t-inv-non-comp} 
Assume that the  generic conditions {\bf (G1)} and {\bf (G3)} hold. 
If 
$\qq$ contains an unbounded inverse component, then the switched
system \r{s-system} is \bad.
\et
\subsection{Robustness}\label{s-robustness}
We say that a property satisfied by  $(X,Y)$ is \underline{robust} if it still holds for small 
perturbations of the pair $(X,Y)$, that is, if it holds for all the elements of 
 a neighborhood of $(X,Y)$ in GAS$(\R^2)\times$GAS$(\R^2)$.
Such notion of robustness is also known as {\it structural stability}, an expression which we prefer to avoid, in order to prevent confusion with the many 
definitions of stability already introduced for \r{s-system}.

Under the generic conditions {\bf (G1)} and {\bf (G2)}, 
one can easily verify that the 
topology of the set $\qq$ does not change 
for small perturbations 
of $X$ and $Y$. Moreover, fixed one 
component 
$\Gamma$ of $\qq$,  the fact that $\Gamma$ is direct or 
inverse is robust. Similarly, if $\Gamma$ is a component of $\qq$, which has not the origin in its 
closure,  the absence of tangency points along $\Gamma$ is robust.
As a consequence, the conditions formulated by the theorems above are robust. More precisely: 

\bt
Under generic assumptions, if any of Theorems \ref{t-mai-paralleli}, \ref{t-mai-tangent}, \ref{t-solo-compatte}, \ref{t-inv-non-comp},  Corollary \ref{c1}, or
Proposition \ref{p-mai-inverse}   
applies to the pair $(X,Y)$, then it applies in a neighborhood of $(X,Y)$ in GAS$(\R^2)\times$GAS$(\R^2)$. 
\et

\section{Proof of Theorem \ref{t-mai-paralleli}}\label{mp}

Assume that $\qq=\{0\}$. We already recalled in Section \ref{s-basic-facts} that GAS and GUAS are two equivalent notions. 
The main step of the proof consists in showing that (\ref{s-system}) is globally attractive.
The uniform stability will be obtained as a byproduct of
the adopted demonstration technique.

Fix $q\in \R^2\setminus\{0\}$. 
We first prove that ${\A(q)}$
is bounded. 
Then we  show that, for every $u(.)$ in $\U$,  
 the only possible accumulation point of $\ga(q,u(.),t)$ is
the origin. 
These two facts imply that $\ga(q,u(.),t)$ 
converges to the origin
as $t$ goes to infinity.

\subsection{Boundedness of ${\A(q)}$}\label{barbecue}

\noindent We distinguish two cases.

\noindent
{\bf \underline{First case:} $\gamma_X(q,.)$ and $\gamma_Y(q,.)$ do not \fintersect}. Then, we
can define a closed, simple, piecewise smooth curve, by
$$
\gamma_{X,Y}(q,t)=\left\{\begin{array}{lr}\!\!\gamma_X(q,\tan(t\pi)) & \mbox{if }\ t\in\left[0,\frac{1}{2}\right],\\[1mm]
\!\!\gamma_Y(q,\tan((1-t)\pi)) & \mbox{if }\ 
t\in\left[\frac{1}{2},1\right],
\end{array}\right.
$$
where $\gamma_X(q,\tan(\pi/2))$ and $\gamma_Y(q,\tan(\pi/2))$ are 
identified with the origin. 
The support of $\gamma_{X,Y}(q,.)$ separates $\R^2$ in two sets, one being
bounded. Let us call $\B(q)$ the interior of the bounded set and 
$\D(q)$ the interior of the unbounded one. 

\bl
\label{goo}
$\A(q)$ is contained in
$\Bbq=\B(q)\cup\gamma_{X,Y}(q,[0,1])$.
\el

\noindent\proof 
Consider the vector field $\dbt$. At the point $q$, it points either inside or outside $\B(q)$. Then, as it becomes clear through a local rectification of $\dbt$, 
the same 
holds true at all points of  $\gamma_{X,Y}(q,[0,1])$ sufficiently close to $q$.
Moreover, 
since the orientation defined by $(X,Y)$ does not vary on 
$\R^2\setminus\{0\}$ and coincides with the
ones induced by $(X,(X+Y)/2)$ and $((X+Y)/2,Y)$, 
then $\dbt$ is pointing  constantly  either inside or outside $\B(q)$, all along 
$\gamma_{X,Y}(q,[0,1])\setminus\{0\}$.

Let us assume that $\dbt$ points inside $\B(q)$. Then $\Bbq$ is 
invariant for the flow of all the vector fields of the type $u\,X+(1-u)\,Y$, with $u\in[0,1]$, that is, it is
invariant for the dynamics of \r{s-system}.
Hence, $\A(q)$ is contained in $\Bbq$.

Assume now, by contradiction, that $\dbt$ points outside $\B(q)$.
The same reasoning as above shows that $\A(q)$ is contained in $\Dbq$. 
Define, for every $t\geq 0$ and every $\tau\in \R$,
$$
\gamma^{X,Y,t}(q,\tau)=\left\{\begin{array}{lr}
\!\!\gamma_X(q,-\tau) & \mbox{if }\ \tau<-t\,, \\
\!\!\gamma_Y(\gamma_X(q,t),\tau+t) & \mbox{if }\ \tau>-t\,.
\end{array}\right.
$$
The support of $\gamma^{X,Y,t}$ is  given by the union of the integral curves of $X$ and $Y$ 
connecting $\gamma_X(q,t)$ and the origin 
(see Figure~\ref{gxyt}). For every $t\geq0$, 
we can identify $\gamma^{X,Y,t}$ with a closed curve passing through the origin. 

\bfi[h]
\begin{center}
\input{inde.pstex_t}
\caption{The curves $\gamma^{X,Y,t}$}\label{gxyt}
\end{center}
\efi

Fix a point $q'$ in $\B(q)$. By hypothesis, no $\gamma^{X,Y,t}$ passes through $q'$.
Notice that the index of
$\gamma^{X,Y,0}$ with respect to $q'$ is equal to one, since the 
support of $\gamma^{X,Y,0}$ coincides with the boundary of $\B(q)$. 
The stability of $Y$ at the origin implies that the index with respect to $q'$ of the curve
$\gamma^{X,Y,t}$ depends continuously of $t$, that is, it is constant on $[0,\infty)$. 
Hence, for every $t\in [0,\infty)$,
\[\max_{\tau\in\R}\|\gamma^{X,Y,t}(q,\tau)\|> \|q'\|>0\,.\]
On the other hand, when $t$ goes to infinity, $\gamma_X(q,t)$ converges to the origin and
\bqn
\sup_{\tau<-t}\|\gamma^{X,Y,t}(q,\tau)\|=\sup_{\tau>t}\|\gamma_X(q,\tau)\|\stackrel{t\rw\infty}{-\!\!\!-\!\!\!-\!\!\!\!\!\longrightarrow} 
0\,.\eqnn
Therefore, there exist $p\in\R^2$, arbitrarily close to the origin, 
such that 
the curve $s\mt \gamma_Y(p,s)$, $s>0$, exits the ball $B_{\|q'\|}$, which contradicts the 
stability of $Y$ at the origin. \EOP

\medskip

\noindent{\bf \underline{Second case:} $\gamma_X(q,.)$ and $\gamma_Y(q,.)$ do \fintersect}. Let $t$ be the first positive time
such that the point $\gamma_Y(q,t)$ is equal to $\gamma_X(q,\tau)$ for some $\tau>0$. Define, for every $s\in [0,\tau+t]$,
$$
\gamma_{X,Y}(q,s)=\left\{\begin{array}{lcr}\!\!\gamma_X(q,s) & \mbox{if }&
s\in[0,\tau]\,,\\[1mm]
\!\!\gamma_Y(q,t+\tau-s) & \mbox{if }&
s\in[\tau,\tau+t]\,.
\end{array}\right.
$$
The curve $\gamma_{X,Y}(q,.)$ is simple and closed, and separates $\R^2$ in two open sets $\B(q)$ and $\D(q)$, $\B(q)$ being
bounded. 

\bl
\label{hoo}
$\A(q)$ is contained in
$\Bbq=\B(q)\cup\gamma_{X,Y}(q,[0,\tau+t])$.
\el
\proof
Assume that
$\dbt$ points inside $\B(q)$ at $q$. Hence, the same is
true for all points of 
$\gamma_{X,Y}(q,[0,\tau+t])$ sufficiently close to $q$. 
Since $\qq=\{0\}$, the property extends to 
the entire curve $\gamma_{X,Y}(q,[0,\tau+t])$, except 
possibly at the point $\gamma_{X,Y}(q,\tau)$. 
%
The same reasoning 
can be applied 
at $\gamma_{X,Y}(q,\tau)$, showing that $\dbt$ points either inside or outside $\B(q)$ at all points
of the type $\gamma_{X,Y}(q,s)$, with $s$ close to $\tau$. The only non-contradictory 
possibility is that
 $\dbt$ points inside $\B(q)$ all along $\gamma_{X,Y}(q,[0,\tau+t])$.
Hence, $\Bbq$ is invariant under the flow of each vector field $u\,X+(1-u)\,Y$, $u\in[0,1]$, so that 
$\A(q)$ is contained in $\Bbq$.

Assume, by contradiction, that $\dbt$ points 
inside $\D(q)$.
The same reasoning as above shows that $\A(q)$ is contained in $\Dbq$.
In particular, the origin belongs to $\Dbq$. 
On the other hand,  the fact that $\dbt$ 
points outside $\B(q)$ all along
$\partial \B(q)$ implies that it has a zero inside $\B(q)$. Which is impossible, unless $\B(q)$ contains the origin. \EOP

\medskip

We proved that, in both cases, the set $\A(q)$ is bounded.
The precise description of $\A(q)$ is given by the following lemma, where the definition of $\B(q)$
depends on whether $\ga_X(q,.)$  and $\ga_Y(q,.)$  \fintersect\ or not.

\bl
\label{gooh}
 $\A(q)=\Bbq\setminus\{0\}$.
\el

\noindent\proof
First notice that the origin does not belong to $\A(q)$, being a steady point for both $X$ 
and $Y$. The inclusion of $\A(q)$ in $\Bbq\setminus\{0\}$ is thus a consequence of Lemma~\ref{goo} 
and Lemma~\ref{hoo}. 

As for the opposite inclusion, notice that $\partial\B(q)\setminus\{0\}$ is, by construction, made of integral curves of $X$ and $Y$ starting from $q$. Therefore, $\partial\B(q)\setminus\{0\}\subset \A(q)$.

Fix now $p\in\B(q)\setminus\{0\}$. We are left to prove that $p\in \A(q)$. 
Define
$$C=\{\gamma_{X}(p,\tau)\,|\;\tau\leq 0\}\,,$$
and let $V$ be a neighborhood  of the origin such that $p\not\in V$. 

Due to the  stability of $X$ and the boundedness of $\B(q)$, there exists $T>0$ such that $\ga_X(\Bbq,T)\subset V$. Since $\ga_X(\gamma_{X}(p,-T),T)=p\not\in V$, then $C$ is not contained in $\Bbq$. Therefore, there exists $\tau<0$ such that $\gamma_{X}(p,\tau)\in\partial \B(q)$. Notice that $\gamma_{X}(p,\tau)$ is different from the origin, since otherwise we would have $p=0$. 
Finally, $\gamma_{X}(p,\tau)\in \A(q)$, which implies that $p=\ga_X(\gamma_{X}(p,\tau),|\tau|)$ belongs to $\A(q)$.  
\EOP

\subsection{Global attractivity}

In the previous section, we showed that the accessible set from every  point is bounded. Hence, 
the global attractivity of (\ref{s-system}) is proved if we ensure that
no admissible curve 
has an accumulation point different from the origin.

Let us show that, for every point $p\neq0$, there exist 
$\eps>0$ and a neighborhood $V_p$ of $p$ such that every admissible curve $t\mapsto \ga(q,u(.),t)$ 
entering $V_p$ 
at time $\tau$ leaves $V_p$ before time $\tau+\eps$ and never comes back to $V_p$ after 
time $\tau+\eps$.


Since $X$ and $Y$ are not parallel at $p$, we can choose a coordinate system $(x,y)$ such that 
$X(p)=(1,-1)$ 
and $Y(p)=(1,1)$. We denote $p=(p_x,p_y)$, $X(x,y)=(X_1(x,y),X_2(x,y))$, $Y(x,y)=(Y_1(x,y),Y_2(x,y))$.
The fields $X$ and $Y$ being continuous, there exists $\al>0$ such that, if $(x,y)\in B_\infty(\al)=
\{(a,b)\;|\; |a-p_x|<\al, \; \; |b-p_y|<\al \}$, then $X_1(x,y)$, $Y_1(x,y)$, $-X_2(x,y)$, and  
$Y_2(x,y)$ are in 
$[{1}/{2},{3}/{2}]$.


Let $p'=(p_x-\frac{\al}{10},p_y)$ and consider $\ga_X(p',.)=(\ga_X^1(p',.),\ga_X^2(p',.))$. Its first 
coordinate $\ga_X^1(p',.)$ is increasing and its derivative takes values in 
$[{1}/{2},{3}/{2}]$. 
The same is true for $-\ga_X^2(p',.)$. Hence $\ga_X(p',.)$ does not leave the set 
$B_\infty(\al)$ before time ${2\al}/{3}$. Since $\ga_X^1(p',2\al/5)$ is larger than
$p_x+\frac{\al}{10}$ and 
$\ga_X^2(p',2\al/5)$ is in $[p_y-\frac{3\al}{10},p_y-\frac{\al}{10}]$, then the curve $\ga_X(p',.)$
intersects the segment $S_p=B_\infty(\al)\cap\{(x,y)|x=p_x+\frac{\al}{10}\}$ 
in a time $\tau_X$ smaller than ${2\al}/{5}$. 

The same occurs for $\ga_Y(p',.)$. Denote by $\tau_Y$ its intersection time with $S_p$.

Choose as $\TT$ the bounded set whose boundary 
is given by the union of
$\ga_X(p',[0,\tau_X])$, $\ga_Y(p',[0,\tau_Y])$, and
the segment $[\ga_X(p',\tau_X),\ga_Y(p',\tau_Y)]=\{\lb\,\ga_X(p',\tau_X)+(1-\lb)\,\ga_Y(p',\tau_Y)\,|\;0\leq \lb\leq 1\}$ (see Figure~\ref{TTT}).

\bfi[h]
\begin{center}
\input{stab2.pstex_t}
\caption{The set $\TT$}\label{TTT}
\end{center}
\efi

The following 
lemma states that $\TT$ satisfies the required properties. 
As a consequence, $p$ cannot be the 
accumulation point of any admissible curve.

\bl \label{pT}
We have the following:
(i) $\TT$ is a neighborhood of $p$; (ii) every admissible curve entering $\TT$ leaves $\TT$ in a time smaller than ${2\al}/{5}$ through  the segment $[\ga_X(p',\tau_X),\ga_Y(p',\tau_Y)]$; (iii) once an admissible curve leaves $\TT$, it enters $\A(p')\setminus \TT$ and never leaves it.
\el

\noindent\proof
The first point follows by the construction of $\TT$. As for ({\it ii}),
notice that all the points of $\TT$ have first coordinate in 
$[p_x-\frac{\al}{10},p_x+\frac{\al}{10}]$. Since the first coordinate 
of $X$ and $Y$ is larger than ${1}/{2}$, then every admissible curve entering $\TT$ 
leaves it in a time smaller than ${2\al}/{5}$. 
Moreover, since along $\ga_X(p',.)$ and $\ga_Y(p',.)$ the admissible velocities of \r{s-system} point inside $\TT$, then an admissible curve can leave $\TT$ only through  the segment $[\ga_X(p',\tau_X),\ga_Y(p',\tau_Y)]$.
Finally, ({\it iii}) follows from the remark that $\A(p')\setminus \TT$ is invariant 
for the dynamics, since the admissible velocities of \r{s-system} 
point inside $\A(p')\setminus \TT$ all along its boundary. \EOP


\subsection{Conclusion of the proof of Theorem \ref{t-mai-paralleli}}
We are left to prove that \r{s-system} is uniformly stable. To this extent, fix $\de>0$.
Since both $X$ and $Y$ are stable at the origin, then
there exists
$\eps>0$ such that every integral curve of $X$ or $Y$ starting in $B_\eps$
is contained in $B_\de$. 
Hence, for every $q\in B_\eps$,
the boundary of $\A(q)$ is contained in  $B_\de$.
Therefore,  $\A(q)$, being bounded, 
is itself contained in $B_\de$. \EOP

\begin{rmk}
The proof of Theorem~\ref{t-mai-paralleli} naturally extends  to the following case: if $V$ is an open and
simply connected subset of $\R^2$, if $X$ and $Y$  point inside $V$ along its 
boundary, and if $\qq\cap V=\{0\}$, 
then 
\r{s-system}
is uniformly asymptotically stable on $V$. 
\end{rmk}

\section{Proof of Theorem \ref{t-mai-tangent}}\label{mt}

The proof follows the main steps as the one of 
Theorem \ref{t-mai-paralleli}. The idea is again to fix a point $q\in\R^2$, to characterize the boundary 
of its accessible set $\A(q)$, to prove that such set is bounded, 
and, finally, to show that no admissible curve has an accumulation point different from the origin.

In order to describe the boundary of ${\A(q)}$, we need some extra construction.
Notice that every component $\Ga$ of $\qq$ separates the plane in two 
parts. 
Since $\Ga$ contains no tangency points, then one of such two regions must be invariant for $X$, and the same argument holds for $Y$ as well.
Necessarily, the invariant region is the one containing the origin, which is attractive both for $X$ and $Y$. In particular, $\Ga$ is direct
and every admissible curve 
crosses $\Ga$ at most once.
Associate with every point $q\in \R^2$ the number $n(q)$ 
of components of $\qq$ that the curve $\ga_X(q,.)$ crosses 
at strictly positive times, 
before converging to the origin (see Figure~\ref{nq}). 
Since the curve $\ga_X(q,(0,\infty))$ is bounded and crosses 
each component of $\qq$ at most  once, then $n(q)$ 
is finite. 
\bfi[h]
\begin{center}
\input{nq.pstex_t}
\capt\label{nq}
\end{center}
\efi


For every $i\leq n(q)$, let us denote by $\Ga_i$ the $i$-th component of $\qq$ crossed by 
$\ga_X(q,.)$. We claim that $\ga_Y(q,.)$ crosses exactly the same components as $\ga_X(q,.)$, 
in the same order. 
Otherwise, as one can easily check, $X$ and $Y$ would not both be GAS at the origin (the reason is that
the components of $\qq$ separate the plane and can be crossed by
an admissible curve at most once).

Let us define two admissible curves, starting from $q$, 
that can be used to characterize the boundary of $\A(q)$, 
in analogy with what has been done in the proof of Theorem \ref{t-mai-paralleli}. 
The first of such curves follows the flow of $X$ until it reaches $\Ga_1$, 
then follows the flow of $Y$ 
until it crosses $\Ga_2$, and so on.
The second one follows alternatively the flows of $X$ and $Y$ in the other way round, 
starting with $Y$ and switching to $X$ as it meets $\Ga_1$. 
Such two curves converge to the origin, since $n(q)$ is finite. 
As in the proof of Theorem~\ref{t-mai-paralleli}, we can distinguish two cases, depending on whether the two curves intersect or not (see Figure~\ref{acc3}). 
\bfi[h]
\begin{center}
\input{acc3.pstex_t}
\capt\label{acc3}
\end{center}
\efi
The arguments of Section~\ref{mp} can be adapted
in order to prove the boundedness of \r{s-system} and the absence of accumulation points different from the origin. The details are left to the reader.\EOP

\section{Proof of Theorem \ref{t-solo-compatte}}\label{sc}

Consider a system of coordinates $(x,y)$ on $\R^2$ which preserves the origin and renders 
$X$ radial outside a ball $B_{R_0}$, 
$R_0>0$. (Such system can be defined using the level sets of a smooth Lyapunov 
function for $X$, see \cite{g2g}.)
Taking possibly a larger $R_0$, we can assume that $X$ and $Y$ are never collinear in $\R^2\setminus B_{R_0}$.

For every $R>0$, 
let
$$\Om_R=\cup_{p\in B_R} \A(p)\,.$$
Our aim is to prove that each $\Om_R$ is bounded. 

Fix $R>R_0+1$.
If $(X,Y)$ is replaced with a pair of vector fields $(\Xt,\Yt)$ which coincides with $(X,Y)$ outside $B_{R_0+1}$, then the set $\Om_R$, constructed as above,
does not change.
The idea is to choose $\Xt$ and $\Yt$ in such a way that they are never parallel outside the origin and still GAS. The boundedness of $\Om_R$ follows then from Theorem \ref{t-mai-paralleli}.

Set
\brs
X_0(x,y)&=&-x\partial_x-y\partial_y\,,\\
Y_0(x,y)&=&y\partial_x-x\partial_y+\lambda X_0\,,\ \ \ \ \ \ \ \ \ \lb>0\,,
\ers
and notice that $X$ and $X_0$ are collinear outside $B_{R_0}$.
Notice, moreover, that, if $\lambda$ is large enough, then the angle between $X_0$ and $Y_0$ is smaller than the minimum of the angles between $X$ and $Y$ in $B_{R_0+1}\setminus B_{R_0}$ (see Figure~\ref{rad2}). Fix such a $\lambda$.

\bfi[h]
\begin{center}
\input{radial2.pstex_t}
\capt\label{rad2}
\end{center}
\efi

The function $Q$ has constant sign on $\R^2\setminus B_{R_0}$. Without lost of generality, we can assume that it is positive. 
Fix a smooth function $\phi:[0,+\infty)\rw[0,1]$ such that $\phi(r)=0$ if $r\leq R_0$ and $\phi(r)=1$ if $r\geq R_0+1$.
Define
\brs
\Xt(x,y)&=&\lp 1-\phi\lp\sqrt{x^2+y^2}\rp\rp X_0(x,y)+\phi\lp\sqrt{x^2+y^2}\rp X(x,y)\,,\\
\Yt(x,y)&=&\lp 1-\phi\lp\sqrt{x^2+y^2}\rp\rp Y_0(x,y)+\phi\lp\sqrt{x^2+y^2}\rp Y(x,y)\,.
\ers

By construction, $(\Xt,\Yt)$ coincides with $(X,Y)$ outside $B_{R_0+1}$ and $det(\Xt,\Yt)$ is strictly positive on $\R^2\setminus\{0\}$. 
We are left to check the global asymptotic stability of $\Yt$,  
the one of $\Xt$ being evident. 
This can be done by using a comparison argument 
between the integral curves of
$Y$ and $\Yt$. 
Indeed, since the angle between $\Xt$ and $\Yt$ is smaller than 
the angle between 
$\Xt$ and $Y$ in $B_{R_0+1}\setminus B_{R_0}$, then 
the integral curve of $\Yt$ starting from a point $q\in\R^2\setminus B_{R_0}$ joins 
$B_{R_0}$ in finite time, 
with a smaller total variation in the angular component that the integral curve of $Y$ starting from the same point $q$.
\EOP

\begin{rmk}
The proof given above applies, without modifications, to the more general case where 
the points at which 
$X$ and $Y$ are 
globally asymptotically stable 
are allowed to be different.
\end{rmk}

\begin{rmk}\label{gus}
The conclusion of Theorem~\ref{t-solo-compatte} would not hold 
under the weaker hypothesis that $X$ and $Y$ are GUS, instead of GAS. 
A counterexample can be given as follows: Let $\varphi:[0,1]\rw\R$ be a smooth function such that $0<\varphi(t)<\pi/2$ for every $t\in(0,1)$ and $\varphi^{(k)}(0)=\varphi^{(k)}(1)=0$ for every $k\geq 0$. Denote by $(r,\theta)$ the radial coordinates on $\R^2$. 
Define, using the radial representation of vectors in $\R^2$, 
\[X(r,\theta)=\left\{\begin{array}{ll}
                                  \lp r,\theta+\frac\pi2+\varphi(r)\rp&\mbox{if }\ r\in[0,1]\,,\\[.8ex]
                                  \lp r,\theta+\frac\pi2-\varphi(r-[r])\rp&\mbox{if \ }r>1\,,
\end{array}
\right.\]
and
\[Y(r,\theta)=\left\{\begin{array}{ll}
                                  \lp r,\theta-\frac\pi2-\varphi(2r)\rp&\mbox{if }\ r\in\left[0,\frac12\right]\,,\\[.8ex]
                                  \lp r,\theta-\frac\pi2+\varphi\lp r+\frac12-\left[r+\frac12\right]\rp\rp&\mbox{if }\ r>\frac12\,,
\end{array}
\right.\]
where $[r]$ denotes the integer part of $r$. 
Then, for every $r\geq 1$, $X(r,\theta)$ and $Y(r,\theta)$ are linearly independent, since
the difference between their angular components is given by
\[0<\pi-\varphi(r-[r])-\varphi\lp r+\frac12-\left[r+\frac12\right]\rp<\pi\,.\]
Hence, $\qq$ is compact. On the other hand, 
the feedback strategy
\[u(t)=\left\{\ba{ll}0&\mbox{if }\ r-[r]\in\left[\frac14,\frac34\right),\\
                     1&\mbox{otherwise}
\ea\right.\]
is such that, for every $p\in\R^2\setminus B_{3/4}$,
$\|\gamma(p,u(.),t)\|$ tends to infinity 
as $t$ tends to $\T(p,u(.))=+\infty$.

Notice that the example can be easily modified in such a way that $\qq$ not only is compact, but actually shrinks to $\{0\}$.
It suffices to take $X(r,\theta)=(r,\theta+\psi_X(r))$ and
$Y(r,\theta)=(r,\theta+\psi_Y(r))$, where the graphs of $\psi_X$
and $\psi_Y$ are as in Figure~\ref{psis}.
\bfi[h]
\begin{center}
\input{psis.pstex_t}
\capt\label{psis}
\end{center}
\efi

\end{rmk}

\section{Proof of Theorem \ref{t-inv-non-comp}}\label{inc}

Let $\Gamma$ be an inverse unbounded component of $\qq$ and assume that {\bf (G1)} and {\bf (G3)} hold.
Due to Lemma~\ref{genio}, $\Gamma$  is a one-dimensional submanifold of $\R^2$, which can be parameterized 
by an injective and smooth map $c:\R\rw \R^2$.

Fix a point $p=(p_x,p_y)=c(\tau)$ on $\Ga$. 
According to the results by Davydov (see \cite[Theorem 2.2
]{dav}),
up to a change of coordinates (which, in particular, sets $p_x=0$), the vector fields $X$ and $Y$ can be represented locally 
by one of the following three normal forms
\begin{enumerate}
\item $X(x,y)=(1,x),\ Y(x,y)=(-1,x)$;\label{passing}
\item $X(x,y)=(1,y-p_y-x^2),\ Y(x,y)=(-1,y-p_y-x^2)$;\label{tu1}
\item $X(x,y)=(-1,x^2-y+p_y),\ Y(x,y)=(1,x^2-y+p_y)$.\label{tu2}
\end{enumerate}

\bfi[h]
\begin{center}
\input{forme-normali-finali.pstex_t}
\capt\label{FN}
\end{center}
\efi
\indent Notice that the type \ref{passing} corresponds to the situation in which $X$ and $Y$ are transversal to $\Ga$ at $p$, while \ref{tu1} and \ref{tu2} are the normal forms 
for the case in which $X$ and $Y$ are tangent to $\Ga$ at $p$.

Recall that $p$ is said to have the {\it small time local transitivity}   property  (STLT, for short) 
if, for every $T>0$ and every neighborhood $V$ of $p$, there exists a neighborhood $W$ of $p$ such that every two points in $W$ are accessible from each other within time $T$ by an admissible trajectory contained in $V$.
It has been proved in  \cite[Theorem 3.1
]{dav} that, under the assumption that the system admits a local representation in normal form, 
$p$ has the STLT property 
if and only if it
is of the type \ref{passing}. 
 In particular, if $p$ is of the type \ref{passing}, then there exist $t(p),T(p)>0$ such that, for every $r,s\in (-t(p),t(p))$ there exists an admissible trajectory which steers $c(\tau+r)$ to $c(\tau+s)$ within time $T(p)$.

Assume now that $p$ is a point of the type \ref{tu1} or \ref{tu2}. 
The curve $\Ga$ stays (locally) on one side of the affine line 
$$p+\mbox{span}(X(p))=\{(x,p_y)\,|\;x\in\R\}\,,$$ 
which is 
the affine tangent space to $\Ga$ at $p$. 
Up to a reversion in the parameterization of $\Ga$, 
we can assume that, for every $t$ in a right neighborhood of $\tau$, $X(c(t))$ points into the locally convex part of the plane bounded by $\Ga$ 
(see Figure~\ref{FN}).
%
%
%
It can be easily verified that the two branches of $\Ga\setminus\{p\}$ are connected by 
integral curves of $X$ and $Y$ arbitrarily close to $p$, in the following sense: 
for every $t>0$ small enough,
there exist $\theta,T>0$  such that, for every $r\in(0,\theta)$, both curves 
$s\mt \ga_X(c(\tau+r),s)$
and
$s\mt \ga_Y(c(\tau-r),s)$
intersect $\Ga$ in a positive time smaller than $T$, and the intersection points are of the type $c(\tau+\rho)$, with $0<|\rho|\leq t$.
We can conclude, using the STLT property at points of $\Ga\setminus\{p\}$ close to $p$,  
 that there exists $t(p)>0$ such that,  for every $\mu\in(0,1)$, 
every two points 
of 
$$\Sigma=\{c(\tau+r)|\;\mu\,t(p)<|r|<t(p)\}$$
can be joined by an admissible trajectory 
of time-length bounded by a uniform $T(p,\mu)>0$.

Therefore, given any pair of points $p_i=c(\tau_i),p_f=c(\tau_f)$ on $\Ga$ of type \ref{passing}, there exists an admissible trajectory going from $p_i$ to $p_f$ of time-length smaller that 
$T(c(\tau_1),\mu_1)+\cdots+T(c(\tau_k),\mu_k)$, where 
$$(\tau_1-t(c(\tau_1)),\tau_1+t(c(\tau_1))),\ldots,(\tau_k-t(c(\tau_k)),\tau_k+t(c(\tau_k)))$$
is a covering of the compact segment of $\R$ bounded by $\tau_i$ and 
$\tau_f$, $\mu_1,\ldots,\mu_k\in (0,1)$ are properly chosen 
 and $T(p,\mu)=T(p)$ if $p$ is of type \ref{passing}.
In particular, system (\ref{s-system}) admits  trajectories going to infinity.\EOP

\begin{rmk}
In the non-generic case the statement of Theorem \ref{t-inv-non-comp} is false. A 
counterexample can be 
found even in the linear case. Indeed, consider the vector fields
\bqn
&&X(q)=A\,q,\ \ \ \ \ \ \ \mbox{ where}~~~ 
A=\left(\ba{cc} -1/20& -1/E\\E&-1/20\ea\right),~~~
E=-\frac{201}{200}-\frac{\sqrt{401}}{200}\,,\nn\\
&&Y(q)=B\,q,\ \ \ \ \ \ \ \mbox{ where}~~~ 
B=\left(\ba{cc} -1/20& -1\\1&-1/20\ea\right).
\eqn
The integral curves of $X$ are ``elliptical spirals'', 
while the integral curves of $Y$ are ``circular spirals''. 
The integral curves of $X$ and $Y$ rotate around the origin in 
opposite sense (since $E<0$).
One can easily check that, in this case, the set $\qq$ is a single 
straight line of equation
\bqn
y=-\frac{20}{\sqrt{401}-1}x\,,
\eqn
and its two components are inverse (see Figure~\ref{ex-l}). 

It can be checked by hand 
that 
the switched system defined by $X$ and $Y$ is GUS, although not GUAS 
(see also \cite{SIAM}, Theorem 2.3, case (CC.3)). 

\bfi[h]
\begin{center}
\input{ex-lin.pstex_t}
\capt\label{ex-l}
\end{center}
\efi

\end{rmk}



\newcommand{\auth}[1]{\textsc{#1}}
\newcommand{\tit}[1]{\textit{#1}}
\newcommand{\jou}[1]{\textrm{#1}}
\newcommand{\bff}[1]{{\bf {#1}}}
\newcommand{\bbf}[1]{{\bf {#1}}}

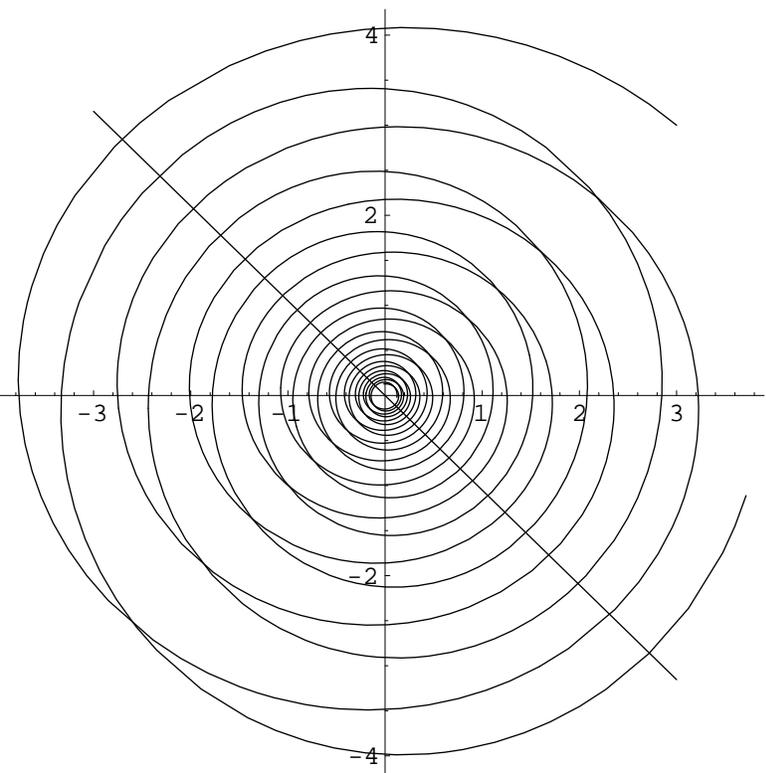
\end{document}